\documentclass[a4paper, 11pt]{article}

\usepackage[english]{babel}
\usepackage[backend=bibtex,
natbib=true,
abbreviate=true,
style=bwl-FU]{biblatex}

\bibliography{/Users/larslau/Dropbox/bib}
\usepackage{amsthm, amsfonts, amssymb, mathrsfs, amsmath}
\usepackage{color}
\usepackage{hyperref} 
\usepackage{mathtools} 
\usepackage{subfig}

\hypersetup{colorlinks=true,%
linkcolor=gray,%
citecolor=gray,%
urlcolor=gray}

\usepackage{tikz}
\usetikzlibrary{arrows,%
                petri,%
                topaths,%
                shapes,%
                snakes,%
                automata,%
                backgrounds}%

\title{\textsc{Differential equations, splines and Gaussian processes}}
\author{Lars Lau Raket}

\newcommand{\ud}{\mathrm{d}}
\newcommand{\by}{\boldsymbol{y}}
\newcommand{\bx}{\boldsymbol{x}}
\newcommand{\bt}{\boldsymbol{t}}
\newcommand{\bs}{\boldsymbol{s}}
\newcommand{\R}{\mathbb{R}}
\newcommand{\I}{\mathbb{I}}

\newcommand{\LL}{\mathscr{L}}
\newcommand{\GG}{\mathcal{G}}
\newcommand{\HH}{\mathscr{H}}
\newcommand{\E}{\mathrm{E}}
\newcommand{\cov}{\mathrm{Cov}}
\newcommand{\btheta}{\boldsymbol{\theta}}

\DeclareMathOperator*{\argmin}{arg\,min}

\theoremstyle{definition}
\newtheorem{example}{Example}[section]

\def\CLOSE{{\boldmath$\circ$}}
\def\markatright#1{\leavevmode\unskip\nobreak\quad\hspace*{\fill}{#1}}
\def\close{\markatright{\CLOSE}}

\def\markatright#1{\leavevmode\unskip\nobreak\quad\hspace*{\fill}{#1}}

\theoremstyle{definition}

\theoremstyle{definition}

\theoremstyle{definition}

\theoremstyle{definition}

\begin{document}
\maketitle
\begin{abstract}
\noindent We explore the connections between Green's functions for certain differential equations, covariance functions for Gaussian processes, and the smoothing splines problem. Conventionally, the smoothing spline problem is considered in a setting of reproducing kernel Hilbert spaces, but here we present a more direct approach. With this approach, some choices that are implicit in the reproducing kernel Hilbert space setting stand out, one example being choice of boundary conditions and more elaborate shape restrictions. 

The paper first explores the Laplace operator and the Poisson equation and studies the corresponding Green's functions under various boundary conditions and constraints. Explicit functional forms are derived in a range of examples. These examples include several novel forms of the Green's function that, to the author's knowledge, have not previously been presented. Next we present a smoothing spline problem where we penalize the integrated squared derivative of the function to be estimated. We then show how the solution can be explicitly computed using the Green's function for the Laplace operator. In the last part of the paper, we explore the connection between Gaussian processes and differential equations, and show how the Laplace operator is related to Brownian processes and how processes that arise due to boundary conditions and shape constraints can be viewed as conditional Gaussian processes. The presented connection between Green's functions for the Laplace operator and covariance functions for Brownian processes allows us to introduce several new novel Brownian processes with specific behaviors. Finally, we consider the connection between Gaussian process priors and smoothing splines.

This paper was originally developed as part of the teaching material for a graduate course in functional data analysis held in 2015 at Department of Mathematical Sciences, University of Copenhagen.
\end{abstract}
\newpage
\begin{quote} 
``Polynomials are wonderful even after they are cut into pieces, but the cutting must be done with care.''\\[-2em]
\flushright{---Isaac Jacob Schoenberg}
\end{quote}

\section{Introduction}
\noindent Splines are piecewise polynomials that are required to have certain smoothness properties. The introductory quote, by one of the early pioneers in the development of splines, encapsulates much of the feeling of working with splines. Splines are truly wonderful and are applicable to a large body of problems, but to understand the details of how they arise and how they are connected to various branches of mathematics and statistics require care.

During the last half century, smoothing splines have become a fundamental tool for estimating continuous functions from discretely observed data. In its typical formulation, the smoothing spline problem involves finding the function $\theta$ that minimizes a discrete Gaussian likelihood term that measures the deviation from observed data plus a penalty term that measures the roughness of the function. The conventional choice of roughness measure is in terms of the integral of its squared second derivative. Let $y_1,\dots,y_m$ denote observed data with associated continuous covariate values $t_1,\dots,t_m$. In this note, we will think of the $t$-variables as representing time and denote them as such, but in practice they can represent any continuous covariate. The classical smoothing spline estimator is defined as the minimizer 
\begin{align}
\hat\theta = \argmin_{\theta} \sum_{i=1}^m (y_i - \theta(t_i))^2 +\lambda \int_T |\theta''(t)|^ 2\,\ud t.\label{eq:smoothing_spline_orig}
\end{align}
It was shown by \cite{holladay1957smoothest} that the function that minimizes the roughness term $\int_T |\theta''(t)|^ 2\,\ud t$ under the constraint of interpolating $y_1,\dots, y_m$ at time points $t_1<\dots < t_m$ is a cubic spline with the so-called natural boundary contions $\theta''(t_1) = \theta''(t_m) = 0$. It easily follows that the solution to \eqref{eq:smoothing_spline_orig} must be a natural cubic smoothing spline, in fact, one can replace the quadratic data term with any reasonable pointwise data term and the solution will still be a natural cubic smoothing spline. The argument is simple: Given any candidate solution $\tilde\theta$, choose $\hat\theta$ as the natural cubic spline that interpolates $\tilde\theta$ in the points $t_1,\dots,t_m$. The pointwise data term will produce the same measure of deviation for $\tilde\theta$ and $\hat\theta$, but the roughness term for the interpolating cubic spline $\hat\theta$ will always be less than or equal to the roughness term for $\tilde\theta$.

Since Holladay's results, much work has been devoted to studying the smoothing spline problem \eqref{eq:smoothing_spline_orig} and its generalizations. The theory and methodology for smoothing splines matured greatly when Grace Wahba and collaborators formulated and solved the problem in the setting of reproducing kernel Hilbert spaces (RKHS) \citep{kimeldorf1971some, Wahba}. While the RKHS approach is powerful, it is not easily accessible in terms of the required level of mathematical skill and the generality of the approach may lead to underappreciation of some of the finer details of the considered problems. The avid reader may for example have noticed how the smoothing spine problem \eqref{eq:smoothing_spline_orig} did not mention the space over which the problem was considered. The classical RKHS approach will simply choose the function space to minimize over as the least restrictive space, but from the perspective of a statistician that choice may not align with the knowledge of the problem at hand. For simplicity, the rest of this paper will explore a slightly simpler smoothing spline problem, namely 
\begin{align}
\hat\theta = \argmin_{\theta} \sum_{i=1}^m (y_i - \theta(t_i))^2 +\lambda \int_T |\theta'(t)|^ 2\,\ud t,\label{eq:smoothing_spline_simple}
\end{align}
where the roughness penalty is replaced by the integrated squared derivative of the function. The solution to this problem will in general be a piecewise linear function with knots in $t_1,\dots,t_m$ as shown by \cite{schoenberg1964spline}, but we show here how we may get complex behavior of the solution by imposing different criteria on the space we are minimizing over. We even show that by imposing sufficiently strict criteria, we may get pathological behavior of the spines, for example such that they become infinitely smooth polynomials. Most of the approaches and connections presented here will seamlessly transfer to the classical smoothing spline problem  \eqref{eq:smoothing_spline_orig} and the general L-spline problem. While requiring some determination, explicit forms of various Green's functions defining the solutions to the classical spline problem can be derived (see e.g. \citeauthor{rytgaard2016statistical}, \citeyear{rytgaard2016statistical}.)

The approach taken here is related to and influenced by a number of previous works. One of them is the work by \cite{Markussen} where the connection between Green's functions and smoothing splines is used to develop a framework for approximate inference for mixed smoothing spline models. Along the way, \citeauthor{Markussen} presents a result \citep[ Theorem 1]{Markussen} that gives an explicit (although quite complex) form for Green's function associated with constant-coefficient linear differential operators under various boundary conditions. This approach was extended to data on multidimensional domains such as $[0,1]^d$ by \cite{RaketMarkussen}. In this paper, the connection between Green's function for various differential operators and covariance functions for known Gaussian processes such as the Brownian sheet, the tied-down Brownian sheet and the Mat\'ern processes are explored in examples. Finally, parts of the work presented here has been extended to the original smoothing spline problem \eqref{eq:smoothing_spline_orig} in \cite{rytgaard2016statistical}. In this work, \citeauthor{rytgaard2016statistical} presents explicit forms for the Green's functions under all combinations of Dirichlet and Neumann boundary conditions along with a number of other extensions involving shape restrictions and alternative likelihood terms.

\section{Green's functions and the Poisson equation}
In this note we will consider the Poisson equation as our driving example. The equation is given by 
\begin{align}
-\partial_t^2f(t) = h(t)\label{eq:poisson}
\end{align}
where $-\partial_t^2$ is the \emph{Laplace operator} $-\partial_t^2f = -f''$, $f\,:\,[0,1]\rightarrow \R$ is twice differentiable in a suitable sense and $h\in L^2([0,1])$. This differential equation is not too difficult to solve, an obvious idea is to write the solutions $f$ as double integrals of $h$
\[
f(t) = -\int_{a_1}^t \int^{s_1}_{a_2} h(s_2)\,\ud s_2\,\ud s_1.
\]
Here we are however not too interested in the specific solutions to the differential equation, but rather in the general problem of inverting the differential operator $-\partial_t^2$. Suppose an inverse operator $\mathscr{G}$ exists, then it should hold that
\[
-\partial_t^2\mathscr{G} h(t) = h(t)
\]
and thus this operator provides a general method for solving the differential equation \eqref{eq:poisson}.

\subsection{Green's functions}
When a differential operator $\LL$ is invertible its inverse $\mathscr{G}$ will be an integral operator with kernel $\GG$,
\begin{align}
\mathscr{G}f(t) = \int_0^1 \GG(s, t) f(s)\,\ud s.\label{eq:green1}
\end{align}
The kernel $\GG $ is referred to as the \emph{Green's function} for $\LL$. Green's functions have two important properties that we will use extensively, but not prove here: firstly, Green's functions are symmetric $\GG (s, t) = \GG (t, s)$; and secondly, the Green's function as a function of the individual coordinates $\GG(s,\,\cdot\,)$ typically obey boundary conditions on the space of function under consideration. In the cases we will consider here, the kernel will be a well-behaved function, but it will sometimes have a more complex structure. In general the Green's function will be a \emph{distribution} or \emph{generalized function}. The most basic distribution is the Dirac delta at $t$, $\delta_t$, for which it holds that  
\begin{align}
\int_0^1 \delta_t(s) f(s)\,\ud s = f(t).\label{eq:dirac}
\end{align}
It is important to note here that the above definition of the Dirac delta function only defines it as the point evaluation functional in combination with an inner product. Thus, it does not always have the standard form of the Dirac delta where one would think of it as a function with $\delta_t(s)=0$ when $s\neq t$ with an infinitely large value at $t$ corresponding to a  point mass. 
\begin{example}\label{ex:indicator}
Using distributions, it is possible to define derivatives of discontinuous functions in a distributional sense. For example, suppose that the constant function $1(t) = 1$ is in the function space we are considering. Then for $a\in [0, 1]$
\[
\int_0^t \delta_a(s)\,\ud s =\boldsymbol{1}_{[a,1]}(t)
\]
where $\boldsymbol{1}_{A}(t)$ is the indicator function that is 1 when $t\in A$ and $0$ otherwise. Thus, in the distributional sense
\[
\partial_t\boldsymbol{1}_{[a,1]}(t) = \delta_a(t).
\]
\close
\end{example}

By applying the differential operator $\LL$ on both sides equation \eqref{eq:green1}, we get the important property of Green's functions that
\[
f(t) = \int_0^1 \LL\GG(s, t) f(s)\,\ud s,
\]
which gives rise to an alternative definition of Green's functions: a Green's function $\GG$ for $\LL$ is a distribution for which
\begin{align}
(\LL\GG(s, \,\cdot\,))(t)=\delta_t(s)\label{eq:green2}
\end{align}
or equivalently
\[
\LL\GG(s, \,\cdot\,)=\delta_s.
\]
Again, we note that this definition is in the distributional sense, thus, depending on the function space we are considering, the representation of $\delta_t$ may vary. 

\begin{example}
From Example~\ref{ex:indicator} we see that $\GG(s, t)=\boldsymbol{1}_{[s,1]}(t)$ is a Green's function for $\partial_t$. 
\close
\end{example}

Since $\LL$ is a linear operator on a function space $\HH$, it is invertible if and only if its kernel is trivial
\[
\ker (\LL) = \{f\in \HH \,:\, \LL f = o\} = \{o\}
\]
where $o(t) = 0$ is the zero function. 

\begin{example}
The kernel for the Laplace operator on the space of twice continuously differentiable functions $\HH=\mathcal{C}^2([0,1])$ are the functions $f$ for which
\[
-\partial_t^2 f(t) = 0\qquad \textrm{for all }t.
\]
Thus, the kernel consists of all affine functions
\[
\ker (-\partial_t^2) =\{f\in \HH \,:\, f(t) = at + b\}
\]
and the Laplace operator is not invertible on $\HH$.
\close
\end{example}

As the previous example showed, we may need to impose extra conditions on the function space $\HH$ in order to identify a unique Green's function for $\LL$. We may however be able to say something about its general form in certain cases.

\begin{example}\label{ex:lin_green}
Any Green's function for the Laplace operator on a subspace $\HH$ of $\mathcal{C}^2([0,1])$ must obey
\[
-\partial_t^2\GG(s, t)=\delta_s(t).
\]
Assuming there exists $f\in\HH$ such that $-\partial_t^2f =1$, we can use Example~\ref{ex:indicator} to write
\[
-\partial_t^2\GG(s, t)=\partial_t\boldsymbol{1}_{[s,1]}(t).
\]
This means means that
\[
\partial_t\GG(s, t)=-\boldsymbol{1}_{[s,1]}(t) + a_1(s),
\]
so that
\begin{align}
\GG(s, t)=-(t-s)\boldsymbol{1}_{[s,1]}(t) + a_1(s)t + a_0(s).\label{eq:lin_green}
\end{align}
Because of the symmetry of $\GG$, $a_1$ and $a_0$ will generally be first-order affine functions. Again we note that this form is only in the distributional sense, and the exact functional form may differ from this. 
\close
\end{example}

\subsection{Boundary conditions}
The typical solution to the problem of non-trivial kernels is to consider linear subspaces of functions subject to boundary conditions. The two most important types of boundary conditions are the \emph{homogeneous Dirichlet boundary conditions}
\[
f(0)=f(1)=0
\]
and the \emph{homogeneous Neumann boundary conditions}
\[
f'(0)=f'(1).
\]
Given a function space, both of these types of boundary conditions will produce linear subspaces of the original space. One can also consider higher-order Neumann boundary conditions
\[
f^{(k)}(0) = f^{(k)}(1) =0, 
\]
and combinations of boundary conditions, both in the form of requiring multiple homogeneous boundary conditions to hold simultaneously, or by specifying linear combinations of boundary conditions that should be zero.

\begin{example}\label{ex:bb_green}
Consider the Laplace operator $-\partial_t^2$ on the subspace $\HH$ of $\mathcal{C}^2([0,1])$ that arise from imposing homogeneous Dirichlet boundary conditions. Since the kernel of the Laplace operator on $C^2([0,1])$ consists of all affine functions, the kernel of $-\partial_t^2$ on $\HH$ is the affine functions for which the boundary conditions hold. Thus $\ker(-\partial_t^2)=\{o\}$ on $\HH$ and the Laplace operator is invertible. 

The Green's function is in a space of distributions associated to $\HH$, which means that $\GG(s, \,\cdot\,)$ must obey the boundary conditions on $\HH$. Since the conditions of Example~\ref{ex:lin_green} hold (see Exercise 1), the Green's function is given by equation \eqref{eq:lin_green}. By the boundary contions we have that
\[
0=\GG(s, 0)=a_0(s)
\]
and
\[
0=\GG(s, 1)=-(1-s) + a_1(s)
\]
Thus,
\begin{align}
\GG(s, t) = -(t-s)\boldsymbol{1}_{[s,1]}(t) + (1-s)t = s\wedge t - st\label{eq:dirichlet_first}
\end{align}
where $\wedge$ denotes the minimum.
\close
\end{example}

\begin{example}\label{ex:neumann0_green}
Consider again the Laplace operator, but this time on the subspace $\HH$ of $f\in\mathcal{C}^2([0,1])$ where $f$ obeys homogeneous Neumann boundary conditions. We see that the kernel of $-\partial_t^2$ consists of all constant functions, and thus we need an extra condition to ensure invertibility. A simple choice is to impose the restriction that $f(0)=0$, which will make the kernel trivial.
\close
\end{example}

\begin{example}\label{ex:bm_green}
We continue to consider the Laplace operator. This time on the subspace $\HH$ of functions $f$ subject to the mixed boundary conditions $f(0)=f'(1)=0$. Clearly the kernel of the Laplace operator is trivial, so a unique Green's function exists. Since there exists $f\in \HH$ such that $-\partial_t^2f(t) = 1$ (see Exercise 1), we can use Example~\ref{ex:lin_green} to give the form of the Green's function. Imposing the boundary conditions on the general form of the Green's function \eqref{eq:lin_green} we get that
\[
0=\GG(s, 0)=a_0(s)
\]
and 
\[
0=\partial_t\GG(s, 1)=-1 + a_1(s).
\]
Thus, the Green's function is 
\begin{align}
\GG(s, t)=-(t-s)\boldsymbol{1}_{[s,1]}(t) + t= s\wedge t.\label{eq:bm_green}
\end{align}

\close
\end{example}


\subsection{Finding Green's functions using Fourier series}
Recall that the \emph{Fourier series} of a function $f\,:\,[0,1]\rightarrow \R$ is given by
\begin{align}
s_f(t) = \frac{1}{2}a_0 + \sum_{i=1}^\infty a_i\cos(2i\pi t) + \sum_{i=1}^\infty b_i\sin(2i\pi t)\label{eq:fourier}
\end{align}
where 
\[
a_i = 2\int_0^1 f(t)\cos(2i\pi t)\,\ud t\qquad\text{and}\qquad b_i = 2\int_0^1 f(t)\sin(2i\pi t)\,\ud t.
\]
One idea for finding Green's functions is to use a Fourier series representation for them, and use the differential operator and boundary conditions to identify the coefficients. We have that 
\begin{align}
-\partial_t^2s_f(t) = 4\pi^2\sum_{i=1}^\infty a_i i^2 \cos(2i\pi t)+4\pi^2\sum_{i=1}^\infty b_i i^2\sin(2i\pi t).\label{eq:diffeq_fourier}
\end{align}
The sine and cosine functions make up an orthogonal basis for $L^2([0,1])$ and we have that
\[
\int_0^1 \cos(2i\pi t)\cos(2j\pi t)\,\ud t = \begin{cases}
1 & \textrm{ if } i = j = 0\\
1/2 &\text{ if } i=j\neq 0\\
0 & \text{ otherwise}
\end{cases},
\]
\[
\int_0^1 \sin(2i\pi t)\sin(2j\pi t)\,\ud t = \begin{cases}
1/2 &\text{ if } i=j\neq 0\\
0 & \text{ otherwise}
\end{cases},
\]
and that the inner product between any sine and cosine function is zero. To identify the coefficients for the solution $s_f$  of the differential equation \eqref{eq:poisson}
\[
-\partial_t^2 s_f(t) = h(t)
\]
we can multiply both sides by respectively $\cos(j\pi t)$ and $\sin(j\pi t)$ and integrate. In the first case, we get
\[
 4 \pi^2 a_j j^2\int_0^1  \cos(2j\pi t)^2\,\ud t = \int_0^1  \cos(2j\pi t)h(t)\,\ud t.
\]
Solving this in $a_j $ gives
\[
a_j =  \int_0^1 \frac{1}{2j^2\pi^2} \cos(j\pi t)h(t)\,\ud t,
\]
and a similar result holds for $b_j$. By rearranging, the solution $s_f$ thus has the form
\begin{align}
s_f(t) = \frac{1}{2}a_0 + \int_0^1  \sum_{i=1}^\infty\frac{1}{2i^2\pi^2} (\cos(2i\pi s)\cos(2i\pi t) + \sin(2i\pi s)\sin(2i\pi t))h(s)\,\ud s.\label{eq:fourier_solution}
\end{align}
Comparing this to the definition of the Green's function \eqref{eq:green1} it is not too hard to see why this may be useful for finding Green's functions. 
\begin{example}[Constant term]
From the definition of $a_0$, we see that $a_0=0$ when
\[
\int_0^1 f(t) \,\ud t =0.
\]
Thus, if we consider the linear subspace $\HH$ of $f\in \mathcal{C}^2([0,1])$ for which this holds, any Green's function has the form
\begin{align}
\GG(s,t)=\sum_{i=1}^\infty\frac{1}{2i^2\pi^2} (\cos(2i\pi s)\cos(2i\pi t) + \sin(2i\pi s)\sin(2i\pi t)).\label{eq:general_greens}
\end{align}
On the other hand, assume that $a_0 \neq 0$.
Using integration by parts twice, the first term can be written on the form
\begin{align*}
\frac{1}{2}a_0&=\int_0^1 f(s)\,\ud s\\
&= [f(s) (2s + c_1)]_0^1 - [f'(s)(t^2 + c_1 s + c_0)]_0^1+ \int_0^1 (s^2 + c_1 s + c_0)f''(t)\,\ud s\\
&= f(1) (2 + c_1) - f(0) c_1 - f'(1)(1 + c_1 + c_0)+f'(0)c_0\\
&\quad- \int_0^1 (s^2 + c_1 s + c_0)h(s)\,\ud s
\end{align*}
where it has been used that $f''(s) = -h(s)$. Thus, if we can make the first four terms disappear, the Green's function has the form
\[
\GG(s, t) =  \sum_{i=1}^\infty\frac{1}{2i^2\pi^2} (\cos(i\pi s)\cos(i\pi t) + \sin(i\pi s)\sin(i\pi t)) - (s^2 + c_1(t) s + c_0(t)).
\]

\close
\end{example}
We note once again that the found forms of the Green's functions are only to be understood in distributional sense, that is, they give the right results  under the integration. But as a function of $t$ the given form of a Green's function will typically not belong to $\HH$. We end this section with an example that demonstrates how to find the closed-form solution of the cosine-series in the general form of the Green's function \eqref{eq:general_greens}.

\begin{example}[Cosine series]
Consider the cosine series
\[
c(s,t)=\sum_{i=1}^\infty \frac{1}{2i^2 \pi^2} \cos(2i \pi s) \cos(2i \pi t).
\]
From the product-to-sum formula, we get that
\begin{align}
\cos(2i \pi s) \cos(2i \pi t) = \frac{1}{2}(\cos(2i \pi (s - t)) + \cos(2i \pi (s+t))).\label{eq:cos_prod_to_sum}
\end{align}
We can thus rewrite the representation
\[
c(s, t) = \sum_{i=1}^\infty \frac{1}{4 i^2 \pi^2} \cos(2i \pi (s - t)) + \sum_{i=1}^\infty \frac{1}{4i^2 \pi^2} \cos(2i \pi (s + t)).
\]
Since the cosine function is even $\cos(2i \pi (s - t)) = \cos(2i \pi |s - t|)$. On the other hand, since the cosine function has period $2\pi$, $\cos(2i \pi (s + t)) = \cos(2i \pi (s + t - \boldsymbol{1}_{[1, 2]}(s + t)))$. 

Using some complex analysis and properties of dilogarithms, one can show that for $s \in [0, 1]$
\begin{align}
\sum_{i=1}^\infty \frac{1}{4i^2 \pi^2} \cos(2i \pi s) = \frac{1}{4}(s(s-1) + 1/6).\label{eq:cos_series}
\end{align}
Putting it all together, we get that
\begin{align}
c(s,t)=\frac{1}{12} + \frac{1}{2}(s^2 + t^2 - s\wedge t) - \frac{1}{2}\boldsymbol{1}_{[1, 2]}(s+t)(s + t - 1)
\end{align}
which can also be written on the form
\[
c(s,t)=\begin{cases}
\frac{1}{12} + \frac{1}{2}(s^2 + t^2 - s\wedge t)& \textrm{for } 0 \le s \le 1 - t\\
\frac{1}{12} + \frac{1}{2}(s(s-1) + t(t-1) + 1 - s\wedge t) & \textrm{for } 1 - t < s \le 1
\end{cases}.
\]

\close
\end{example}

\subsection{Fourier bases and linear subspaces}
The approach for finding Green's functions for the Poisson equation presented in the previous section is very general. To find the Green's function $\GG$ for $-\partial_t^2$ on some subspace $\HH$, we may formulate the constraints of the subspace in terms of the Fourier basis functions.
\begin{example}\label{ex:dirichlet2}
Consider the space $\HH$ of $f\in \mathcal{C}^2([0,1])$ subject to homogeneous Dirichlet boundary conditions. Imposing these restrictions on the representation \eqref{eq:fourier} gives that 
\[
0=s_f(0) = \frac12a_0 + \sum_{i=1}^\infty a_i
\]
and
\[
0=s_f(1) = \frac12a_0 + \sum_{i=1}^\infty a_i.
\]
Thus, the representation of $f$ becomes
\[
s_f(t) =\sum_{i=1}^\infty a_i(\cos(2i\pi t) - 1) + \sum_{i=1}^\infty b_i\sin(2i\pi t),
\]
and the Green's function can be written as
\[
\GG(s, t) = \sum_{i=1}^\infty \frac{1}{2 i^2 \pi^2}((\cos(2i\pi s) - 1)(\cos(2i\pi t) - 1)+ \sin(2i\pi s)\sin(2i\pi t)).
\]
From Example~\ref{ex:bb_green} we know that this series must equal $s\wedge t - st$, but we can also show that explicitly. The product-to-sum formulas for the sine functions give that 
\begin{align}
\sin(2i \pi s) \sin(i \pi t) = \frac{1}{2}(\cos(2i \pi (s - t)) - \cos(2i \pi (s+t))).\label{eq:sin_prod_to_sum}
\end{align}
The cosine term $((\cos(2i\pi s) - 1)(\cos(2i\pi t) - 1)$ equals
\[
\cos(2i\pi s) \cos(2i\pi t) + 1 - \cos(2i\pi s) - \cos(2i\pi t),
\]
so using the cosine product-to-sum formula~\eqref{eq:cos_prod_to_sum}, we get that the Green's function has the form
\[
\GG(s, t) = \sum_{i=1}^\infty \frac{1}{2 i^2 \pi^2}(\cos(2i \pi (s - t)) + 1 - \cos(2i\pi s) - \cos(2i\pi t)).
\]
which, using the cosine-series result \eqref{eq:cos_series} gives the explicit form
\[
\GG(s, t) =  \frac{1}{2}\bigg((|s-t|(|s - t|-1) + 1/6)+\frac{1}{6}-(s(s-1) + 1/6)- (t(t-1) + 1/6)\bigg).
\]
which can easily be manipulated to give the form
\[
\GG(s, t) = s\wedge t - st.
\]

\end{example}

\begin{example}[Zero coefficients]
Consider (the closure of) a linear subspace $\HH$ of $C^{2k}$ spanned by Fourier basis functions subject to the constraints
\[
a_{i_1}=0,\dots, a_{i_{\ell}}=0.
\]
Suppose that $-\partial_t^2$ is invertible on on $\HH$.  From \eqref{eq:diffeq_fourier}, we see that any admissible function $h$ on the left-hand side of the differential equation~\eqref{eq:poisson} will also have zero coefficients for the corresponding basis functions. Thus, due to the orthogonality all the corresponding terms will disappear from the series inside the integral in \eqref{eq:fourier_solution}.
\close
\end{example}

\begin{example}[Linear combinations of coefficients]
Consider (the closure of) a subspace $\HH$ of $C^{2k}$ spanned by the Fourier basis functions subject to a linear constraint on the coefficients
\[
\sum_{i=0}^\infty c_ia_i = 0.
\]
where $c_i \in \R$. Suppose that $-\partial_t^2$ is invertible on on $\HH$, and for simplicity assume that $c_0=0$, so we have no constraint on $a_0$. If we consider the general form of the Green's function \eqref{eq:general_greens}, we see that it will generally not be in $\HH$, since this would require that 
\[
\sum_{i=0}^\infty c_i\frac{1}{i^2}=0. 
\]
To correct for this, consider the compensating term
\[
\GG_c(s,t)=\sum_{i=0}^\infty c_i\frac{1}{2i^2\pi^2}\cos(2i\pi s) \cos(2i\pi t).
\]
\[
-\partial_t^2 \GG_c(s,t)=\sum_{i=0}^\infty c_i\cos(2i\pi s) \cos(2i\pi t).
\]
Using the Fourier series representation $s_f$ of a function $f\in \HH$, we see that 
\[
\int_0^1 -\partial_t^2\GG_c(s,t) s_f(s)\,\ud s = \frac{1}{2}\sum_{i=1}^\infty c_i a_i=0.
\]
Thus, we can subtract $\GG_c$ from the general form of the Green's function without it having any effect in the distributional sense.
\close
\end{example}

\subsection{Examples of Green's functions}
\begin{example}[Balanced, periodic functions]
Consider the subspace $\HH$ of functions $f\in \mathcal{C}^2([0,1])$ for which $\int_0^1 f(t) \,\ud t = 0$. The kernel of the Laplace operator on $\HH$ is
\[
\ker(-\partial_t^2)=\{f\in \HH : f(t) = a (t - 1/2)\},
\] 
so we need one additional constraint. The least restricting constraint is that $f(0)=f(1)$, so we restrict ourselves to those functions. The integral constraint is equivalent to $a_0=0$ in the Fourier expansion \eqref{eq:fourier}, and all the Fourier basis functions obey the constraint of periodicity (that the endpoints are the same).
Thus, the Green's function is on the form
\begin{align}
\GG(s, t) = \sum_{i=1}^\infty \frac{1}{2i^2 \pi^2} (\cos(2i \pi s) \cos(2i \pi t) +\sin(2i \pi s) \sin(2i \pi t)),\label{eq:balanced_green}
\end{align}
which we can rewrite using the sine and cosine product-to-sum formulas as we did in Example~\ref{ex:dirichlet2}. This gives the the representation
\[
\GG(s, t) = \sum_{i=1}^\infty \frac{1}{2i^2 \pi^2} \cos(2i \pi (s - t)).
\]
Using the cosine series \eqref{eq:cos_series} we arrive at the closed form solution
\[
\GG(s, t) = \frac12|s-t|^2 - \frac12 |s-t|+\frac{1}{12}.
\]
We note that this Green's function is equivalent to the covariance function of the periodic Brownian motion constructed by \citeauthor{mumford2010pattern} (\citeyear{mumford2010pattern}; Section 6.6.1). 
\close
\end{example}

\begin{example}[Odd functions]
Consider the space of functions $\HH\subset \mathcal{C}^2([0,1])$ for which $f(t) = -f(1 - t)$ for $t\in [0, 1/2]$. Again we see that the kernel for the Laplace operator consists of the lines that are zero at $1/2$. Under the oddness constraint, the additional constraints $f(0)=f(1)$, $f(0) = 0$ and $f(1) = 0$ are all equivalent. Let $\HH$ be the restriction to these functions. 

The oddness implies that $\int_0^1 f(t) \,\ud t = 0$, so the Green's function must have a similar form to the previously found Green's function \eqref{eq:balanced_green}. The terms that are left out are exactly all the cosine functions, since they are even. Using the product-to-sum formula for sine, the Green's function is found to be
\[
\GG(s, t) = \sum_{i=1}^\infty \frac{1}{4i^2 \pi^2}(\cos(2i \pi (s - t)) - \cos(2i \pi (s+t))).
\]
We see that this must equal
\begin{align*}
\GG(s, t) &=\frac{1}{4} \bigg(|s-t|(|s-t| - 1) - (s+t- \boldsymbol{1}_{[1, 2]}(s + t))(s+t - \boldsymbol{1}_{[1, 2]}(s + t) - 1)\bigg)\\
&= \begin{cases}
\frac{1}{2} s\wedge t- st & \text{for } s + t \le 1\\
s\wedge t + \frac12 s\vee t- st  - \frac12& \text{for } s + t > 1
\end{cases}
\end{align*}
\close
\end{example}

\begin{example}[Mixed boundary conditions]
In this example we will consider the Laplace operator $-\partial_t^2$ on function spaces $\HH$ consisting of functions $f\in \mathcal{C}^2([0,1])$ satisfying the mixed boundary conditions $f(0)=f'(1)=0$ and possibly other conditions. In the case that no other conditions are available, we have already seen that the Green's function is
\[
\GG(s, t) = s\wedge t. 
\]
If we add the constraint that 
\[
\int_0^1 f(t)\,\ud t = 0,
\]
we see that the previous Green's function no longer works, since
\[
\int_0^1\GG(s,t)\,\ud t = -\frac{s^2-2s}{2}.
\]
On the other hand, we see that any constant multiplied to $f$ will integrate to zero, and thus, looking at the definition of the Green's function it seems natural to add a quadratic polynomial to the Green's function (which will become a constant after applying the Laplace operator). In particular, we choose the polynomial $p_s(t)=a_2(s) t^2 + a_1(s)t + a_0(s)$ that obey the boundary conditions and integrate to $\frac{s^2-2s}{2}$. We easily see that $a_0(s) = 0$, and from the the boundary condition at 1,
\[
-2a_2(s) = a_1(s).
\]
We now simply need to determine $a_2(s)$ such that
\[
\frac{s^2-2s}{2} = \int_0^1 p_s(t)\,\ud t = \int_0^1 a_2(s) t^2 - 2a_2(s)t \,\ud t  = -\frac{2}{3}a_2(s).
\]
We see that, $a_2(s) = - \frac{3}{4}(s^2 - 2 s)$, and thus the Green's function becomes
\[
\GG(s,t)=(s\wedge t)  - \frac{3}{4}(s^2 - 2 s)(t^2 - 2t).
\]
\end{example}

\begin{example}[Dirichlet boundary]
In this example we will consider the Laplace operator $-\partial_t^2$ on function spaces $\HH$ consisting of functions $f\in \mathcal{C}^2([0,1])$ satisfying homogeneous Dirichlet boundary conditions $f(0)=f(1)=0$ and possibly other conditions. In the simplest case we saw that the Green's function was
\begin{align}
\GG(s, t) = s\wedge t - st. \label{eq:dirichlet_green}
\end{align}
Consider the additional constraint
\[
\int_0^1 f(t)\,\ud t = 0.
\]
As in the previous example, it seems natural to substract a polynomial that is second order in $s$ and $t$ from the Green's function \eqref{eq:dirichlet_green} to find a new Green's function. We see that
\[
\int_0^1 (s\wedge t - st)\,\ud t = \frac{1}{2}(1-s)s.
\]
From the boundary conditions, the polynomial must be of the form 
\[
p_s(t) =a(s) t^2 - a(s)t.
\]
We have that
\[
\int_0^1 a(s) (t - 1) t\,\ud t = -\frac{1}{6}a(s)
\]
so $a(s)=3(1-s)s$. The Green's function becomes
\[
\GG(s, t) = s\wedge t - st - 3(1-s)s (1-t)t.
\]
\close
\end{example}

\begin{example}[Second order polynomials]\label{ex:poly_motion}
Consider the constraint that $f'''(t) = 0$. The subspace of such functions are the second order polynomials.  Now let $\HH$ be the space of second order polynomials subject to the mixed Boundary conditions $f(0)=f'(1)=0$. This means that any function $f\in \HH$ has the form
\[
f(t) = a(t - 2)t ,\qquad a\in \R.
\]
The Green's function should obey that 
\[
\int_0^1 -\partial_t^2\GG(s, t) f(t)\,\ud t = f(s),
\]
but, using integration by parts twice, we get that the above equation is equivalent to
\[
-2a\int_0^1\GG(s, t)\,\ud t = a(s - 2)s.
\]
Thus,
\[
\int_0^1 \GG(s, t) \,\ud t =-\frac{1}{2}(s - 2)s,
\]
meaning that the Green's function is
\[
\GG(s, t) = \frac{3}{4}(s-2)s(t-2)t.
\]
\close
\end{example}

\begin{example}[Second order polynomial bridges]
Consider the same setup as in Example~\ref{ex:poly_motion}, but with the subspace constraint that $\HH$ is the space of second order polynomials subject to homogeneous Dirichlet boundary conditions. This means that any function $f\in \HH$ has the form
\[
f(t) = a(t -1)t,\qquad a\in \R.
\]
The Green's function should obey that 
\[
\int_0^1 -\partial_t^2\GG(s, t) (a(t -1)t)\,\ud t = a(s -1)s,
\]
but we have that
\[
\int_0^1 \partial_t^2\GG(s, t) t \,\ud t =0
\]
and
\[
\int_0^1 \partial_t^2\GG(s, t) t^2 \,\ud t = 2\int_0^1\GG(s, t) \,\ud t.
\]
Thus it must hold that
\[
\int_0^1 \GG(s, t) \,\ud t = -\frac{1}{2}(s -1)s.
\]
Thus the Green's function has the form
\[
\GG(s, t) = 3 s (1-s)t(1-t).
\]
\close
\end{example}

\section{Smoothing splines}
Let pairs of observation times and observations $(t_1, y_1),\dots, (t_m, y_m)$ be given. Assume that the observations are generated according to the statistical model
\begin{align}
y_i = \theta(t_i) + \varepsilon_i\label{eq:model1}
\end{align}
where the $\varepsilon_i$s are independent  $\mathcal{N}(0, \sigma^2)$-variables, and $\theta\in \HH$ where $\HH$ is some space of functions from the unit interval into the reals. The negative log likelihood is 
\begin{align}
\ell(\theta, \sigma^2) = \frac{m}{2}\log \sigma^2 + \frac{1}{2\sigma^2}\sum_{i=1}^m (y_i - \theta(t_i))^2,\label{eq:like}
\end{align}
but if $\HH$ is a flexible space of functions, such as $\mathcal{C}^2([0,1])$ there will be infinitely many functions that matches the observations exactly, so the maximum likelihood estimator is not unique. Furthermore, since we have noise on our observations, we do not want our estimator of $\theta$ to be an interpolating function, but rather a function that matches the trend of the observations while filtering out some of the noise. One way of achieving this is to add a term to the likelihood~\eqref{eq:like} that penalizes roughness. Here we will consider the integrated square derivative magnitude as our choice of roughness measure 
\[
P(\theta)=\int_0^1 \|\theta'(t)\|^2\,\ud t.
\]
We can define the penalized log likelihood functional in $\theta$ as
\begin{align}
\ell_{\lambda}(\theta)= \sum_{i=1}^m (y_i - \theta(t_i))^2 + \lambda P(\theta)\label{eq:penloglike}
\end{align}
where $\lambda > 0$ is a parameter that controls the weighting between data fidelity and regularity. We will call the minimizer of this functional 
\[
\hat\theta = \argmin_{\theta\in \HH} \ell_\lambda(\theta)
\]
a smoothing spline (although it will typically not be smooth). In the following we will go through the necessary theory of identifying the smoothing spline. 

\subsection{Calculus of variations}
Minimization of functionals is the topic of calculus of variations. The functional derivative $\delta F/ \delta f$ of a functional $F\,:\, f\mapsto F(f)$ in the present setting is defined by the equation
\begin{align}
\int_0^1\frac{\delta F}{\delta f}(t) h(t)\,\ud t = \frac{\ud}{\ud \epsilon} F(f + \epsilon h)\Big|_{\epsilon = 0}\label{eq:fderiv}
\end{align}
where $h$ is an arbitrary function living in the same function space as $f$. As with conventional derivatives, a minimizer is found by solving for the functional derivative equal to the zero function.

\begin{example}
Let us find the functional derivative of the roughness measure $P$. Taking the right-hand side of \eqref{eq:fderiv}, we see that
\begin{align*}
\frac{\ud}{\ud \epsilon }\int_0^1 \|(\theta+\epsilon h) '(t)\|^2\,\ud t\Big|_{\epsilon = 0} &=\frac{\ud}{\ud \epsilon }\int_0^1 (\theta'(t) + \epsilon h'(t))^2\,\ud t\Big|_{\epsilon = 0}\\
&=2\int_0^1 \theta'(t) h'(t)\,\ud t.
\end{align*}
Using integration by parts, this integral can be rewritten
\begin{align}
2\int_0^1 \theta'(t) h'(t)\,\ud t =2(\theta'(1) h(1) - \theta'(0) h(0)) -2\int_0^1 \theta''(t) h(t)\,\ud t .\label{eq:roughness_derivative}
\end{align}
Thus, if we for example consider a space of functions satisfying homogeneous Dirichlet boundary conditions $f(0)=f(1)=0$, the functional derivative is $\delta P/ \delta \theta = 2\partial_t^2\theta$. To minimize the roughness measure we thus need a function $\theta$ such that $\partial_t^2\theta(t) =0$, which in the present case amounts to $\theta(t) = 0$.
\close
\end{example}

In the previous example we saw that minimization of functionals may be related to differential equations. There are a number of interesting links between splines, Gaussian processes and stochastic differential equations. In the following we will find the penalized maximum likelihood estimator for $\theta$. First, we need to identify the functional derivative of the data term in the penalized log likelihood \eqref{eq:penloglike}. Computing the right-hand side of \eqref{eq:fderiv} we get
\[
 -2\sum_{i=1}^m (y_i - \theta(t_i)) h(t_i).
\]
As a consequence, the functional derivative must be of the form
\[
 -2\sum_{i=1}^m (y_i - \theta(t))\delta_{t_i}(t).
\]
If we choose boundary conditions such that the boundary terms in \eqref{eq:roughness_derivative} disappear, we can write the functional derivative of $\ell_\lambda$ as 
\[
 -2\sum_{i=1}^m (y_i - \theta(t))\delta_{t_i}(t)-2\lambda \theta''(t).
\]
Suppose the Laplace operator $-\partial_t^2$ is invertible on $\HH$ with Green's function $\GG$. Then we can write $\delta_{t_i} = -\partial_t^2 \GG(t_i, t)$. Solving the for the functional derivative equal to the zero function then leads to 
\[
 \frac{1}{\lambda}\sum_{i=1}^m (y_i - \theta(t))(-\partial_t^2 \GG(t_i, t))=-\partial_t^2\theta(t).
\]
Consider a solution on the form
\[
\theta(t) = \sum_{i=1}^m c_i \GG(t_i, t),
\]
insertion in the differential equation gives 
\[
 \frac{1}{\lambda}\sum_{i=1}^m \bigg(y_i - \sum_{j=1}^m c_j \GG(t_j, t)\bigg)\delta_{t_i}(t)=\sum_{i=1}^m c_i \delta_{t_i}(t),
\]
but since $\delta_{t_i}$ is the point evaluation in $t_i$ we can substitute $\GG(t_j, t)$ with $\GG(t_j, t_i)$. The coefficients $c_i$ can thus be found by solving a linear system of equations
\begin{align}
\sum_{j=1}^m c_j (\GG(t_j, t_i)+ \lambda\boldsymbol{1}_{\{i\}}(j))  = y_i,\qquad i = 1,\dots, m.\label{eq:coef_system}
\end{align}
Define the matrix of the Green's function evaluated at all combinations of observation points
\[
G = \begin{pmatrix}
\GG(t_1, t_1) & \cdots & \GG(t_1, t_m)\\
\vdots & \ddots & \vdots\\
\GG(t_m, t_1) & \cdots & \GG(t_m, t_m)
\end{pmatrix}.
\]
The linear system \eqref{eq:coef_system} can be written on the form
\[
(G - \lambda \I_m)\boldsymbol{c} = \by
\]
where $\I_m$ is the $m\times m$ identity matrix and bold indicates vectorization. Thus, the coefficients are given by
\[
\boldsymbol{c}= (G + \lambda \I_m)^{-1}\by
\]
Finally, define the vector valued function
\[
G(t) = \begin{pmatrix}
\GG(t_1, t)\\
\vdots\\
\GG(t_m, t) 
\end{pmatrix}.
\]
Putting everything together, we have shown that
\begin{align}
\hat\theta (t) = G(t) (G + \lambda \I_m)^{-1}\by.\label{eq:smoothing_spline}
\end{align}
In particular, if we want to evaluate $\hat\theta$ in all observation points, we can compute
\[
\begin{pmatrix}
\hat\theta (t_1)\\
\vdots\\
\hat\theta (t_m)\\
\end{pmatrix} = G (G + \lambda \I_m)^{-1}\by.
\]

\section{Gaussian processes and differential equations}
A Gaussian process $x$ is a continuous-time stochastic process for which all finite-dimensional distributions follow multivariate normal distributions. In this note, we will largely ignore the details of the underlying probability space and think of $x$ as a random variable on $L^2([0, 1])$ such that the individual sample paths $x\,:\,[0,1]\rightarrow \R$ are considered random functions.

\subsection{Brownian motion}
The most important of all stochastic processes is the Brownian motion. We will leave its construction to other courses, but give one definition of it. The Brownian motion $x$ is a stochastic processes for which $x(0)=0$ almost surely, and the increments 
\[
x(t_1)-x(0),\quad x(t_2) - x(t_1),\quad\dots\quad x(t_n) - x(t_{n-1})
\]
are independent, zero-mean normally distributed with variances $t_{i+1} - t_i$ for all choices of $n$ and $0 < t_1 <\dots < t_n$.

From this definition, we can derive some important properties of Brownian motion. Clearly $x$ is a zero-mean process, that is
\[
\E[x(t)] = 0\,\qquad \textrm{for all }t\in[0, 1].
\]
We can also find the covariance function. Let $s \le t$
\[
\cov[x(s), x(t)]=\E[x(s) x(t)]-\E[x(s)]\E[x(t)] = \E[x(s) x(t)]
\]
Now, write $x(t) = x(t) - x(s) + x(s)$. Then
\[
\cov[x(s), x(t)]=\E[x(s) (x(t) - x(s))] + \E[x(s)^2] =s.
\]
We thus conclude that for general $s, t\in [0, 1]$
\[
\cov[x(s), x(t)]=s\wedge t,
\]
and note that we found this expression as the Green's function for the Laplace operator in equation~\eqref{eq:bm_green}. 
\begin{example}[The Brownian bridge]
Let $x$ be a Brownian motion. Let $0\le s \le t \le 1$. Since $x$ is a Gaussian process, the distribution of $(x(s), x(t), x(1))$ is a multivariate normal distribution with mean $
\mu = (0, 0, 0)$ and covariance matrix $\Sigma$ given by evaluating the covariance function at all combinations og the points $s$, $t$ and $1$. That is,
\[
\Sigma = \begin{pmatrix}
s & s & s\\
s & t & t \\
s & t & 1
\end{pmatrix}.
\]
From the theory of the multivariate normal distribution, we know that if $(\bx_1, \bx_2)$ is joint normal with mean $(\boldsymbol{\mu}_1, \boldsymbol{\mu}_2)$ and covariance matrix
\[
\begin{pmatrix}
\Sigma_{11} & \Sigma_{12}\\
\Sigma_{21} & \Sigma_{22}
\end{pmatrix},
\]
then the conditional distribution of $\bx_1 | \bx_2 = \boldsymbol{a}$ is again multivariate normal with mean 
\[
\boldsymbol{\mu}_1 + \Sigma_{12}\Sigma_{22}^{-1}(\boldsymbol{a} - \boldsymbol{\mu}_2)
\]
and covariance matrix
\[
\Sigma_{11} - \Sigma_{12} \Sigma_{22}^{-1} \Sigma_{21}.
\]
Now consider the condition $x(1) = 0$. We can use the above formulas with $\bx_1 = (x(s), x(t))$ and $\bx_2 = x(1)$ to get that the conditional distribution of $(x(s), x(t)) | x(1)=0$ is a zero-mean normal distribution with covariance matrix
\[
\Sigma = \begin{pmatrix}
s - s^2 & s - st\\
s - st & t - t^2
\end{pmatrix}.
\]
This process is called the Brownian bridge and we see that its covariance function must be $s\wedge t - st$. Again we note the connection to the Green's function for the Laplace operator found in Example~\ref{ex:bb_green}.
\close
\end{example}

\subsection{Gaussian priors for curves}
In Bayesian statistics, one can choose prior probability distributions for parameters in the model. The Bayesian  equivalent to maximum likelihood estimation is maximum a posteriori estimation, that is, we maximize the posterior probability given the data and the prior. Assume that we are interested in the parameter $\theta$ and have data $\by$. The posterior probability (distribution) is given by Bayes' theorem
\[
p(\theta\,|\, \by) = \frac{p(\by\,|\,\theta)p(\theta)}{p(\by)}
\] 
where $p(\by\,|\, \theta)$ is the likelihood, $p(\theta)$ is the prior and $p(\by)$ is known as the evidence. The evidence may be very hard to compute because we need to integrate out $\theta$ in the joint distribution, but typically we can avoid computing this term altogether, since it is just a normalization constant that is independent of $\theta$. The maximum a posteriori estimate $\hat\theta_{\mathrm{MAP}}$ can be found by minimizing
\[
-\log p(\theta\,|\, \by) \propto -\log p(\by\,|\,\theta)-\log p(\theta).
\] 
Comparing this to the penalized log likelihood \eqref{eq:penloglike}, we see that both methods minimize the sum of the negative log likelihood and an extra term which in both cases acts as a regularization term. 

Consider again the statistical model \eqref{eq:model1} where we observe a function $\theta \in \HH$ with independent, identically distributed Gaussian noise. We may choose a zero-mean Gaussian process with covariance function $\sigma^2 \tau^2 \GG(\,\cdot\,,\,\cdot\,)$ as our prior for $\theta$. This means that any finite evaluation of $\theta$ in the vector of points $\bs$, $\btheta=(\theta(s_1), \dots, \theta(s_\ell))$ follows a multivariate normal distribution of the form
\[
p(\btheta) = \frac{1}{\sqrt{2^\ell\pi^{\ell} \det G(\bs, \bs)}} \exp\bigg(-\frac{1}{2\sigma^2\tau^2}\btheta^\top G(\bs, \bs)^{-1}\btheta\bigg)
\]
where 
\[
G(\bs, \bs)=\begin{pmatrix}
\GG(s_1, s_1) & \dots & \GG(s_1, s_\ell)\\
\vdots & \ddots & \vdots\\
\GG(s_\ell, s_1) & \dots &\GG(s_\ell, s_\ell) 
\end{pmatrix}.
\]
The joint distribution of $\by$ and $\btheta$ is zero-mean multivariate normal with covariance matrix
\[
\begin{pmatrix}
\sigma^2\tau^2 G(\bt, \bt) + \sigma^2\I_{m\times m} & \sigma^2\tau^2 G(\bt, \bs)\\
\sigma^2\tau^2 G(\bs, \bt) & \sigma^2\tau^2 G(\bs, \bs)
\end{pmatrix}
\]
where $\bt=(t_1, \dots, t_m)$ is the vector of observation points. The maximum posterior estimate of $\theta$ in the points $\bs$ is the conditional expectation of $\btheta$ given $\by$, which is given by
\begin{align}
\hat\btheta_{\mathrm{MAP}} =G(\bs, \bt)(G(\bt, \bt) + \tau^{-2} \I_{m\times m})^{-1} \by.\label{eq:MAP}
\end{align}
Comparing with the previously found smoothing spline \eqref{eq:smoothing_spline}, we see that the maximum a posteriori estimate will be identical if we choose $\GG$ as the Green's function for the Laplace operator on $\HH$ and set $\lambda = \tau^{-2}$.

\subsection{Conditional Gaussian processes}
In this section we will consider a zero-mean Gaussian process $x$ with covariance function $\GG$.

\begin{example}[Neumann boundary conditions]\label{ex:cond_neumann}
Consider the observations of $x$ at the points $s$, $t$, $1-\epsilon$ and $1$, ordered increasingly. We see that $(x(s), x(t), (x(1) - x(1-\epsilon))/\epsilon)$ follow a zero-mean normal distribution with covariance
\[
\begin{pmatrix}
\scriptstyle\GG(s, s) &\scriptstyle \GG(s, t) & \frac{\GG(s, 1) - \GG(s, 1-\epsilon)}{\epsilon}\\[0.2em]
\scriptstyle\GG(t, s) &\scriptstyle \GG(t, t) & \frac{\GG(t, 1) - \GG(t, 1-\epsilon)}{\epsilon}\\[0.2em]
 \frac{\GG(s, 1) - \GG(s, 1-\epsilon)}{\epsilon} & \frac{\GG(t, 1) - \GG(t, 1-\epsilon)}{\epsilon} & \frac{\GG(1, 1) + \GG(1 - \epsilon, 1 - \epsilon) - 2\GG(1, 1 - \epsilon)}{\epsilon^2}\\
\end{pmatrix}.
\]
Thus, the conditional distribution of $(x_s, x_t)$ given that $(x(1) - x(1-\epsilon))/\epsilon = 0$ is zero-mean normal with covariance
\[
\begin{pmatrix}
\scriptstyle\GG(s, s) &\scriptstyle \GG(s, t)\\[0.2em]
\scriptstyle\GG(t, s) &\scriptstyle \GG(t, t)
\end{pmatrix} - 
  \frac{\scriptstyle\epsilon^2} {\scriptstyle\GG(1, 1) + \GG(1 - \epsilon, 1 - \epsilon) - 2\GG(1, 1 - \epsilon)}
 \begin{pmatrix}
\scriptstyle\GG_\epsilon' (s, 1)^2 & 
\scriptstyle\GG_\epsilon' (s, 1)\GG_\epsilon' (t, 1)\\[0.2em]
\scriptstyle\GG_\epsilon' (s, 1)\GG_\epsilon' (t, 1) & 
\scriptstyle\GG_\epsilon' (t, 1)^2
\end{pmatrix}
\]
where $\GG_\epsilon' (s, 1) =\GG_\epsilon' (1, s) = \frac{\GG(s, 1) - \GG(s, 1-\epsilon)}{\epsilon}$.
\close
\end{example}

\begin{example}\label{ex:bm_cond_neumann}
Let $x$ be a Brownian motion. The covariance function is $\GG(s, t) = s\wedge t$. The second term in the formula for the conditional variance from Example~\ref{ex:cond_neumann} is zero, since $\GG_\epsilon' (s, 1)= 0$ for all $s \le 1 - \epsilon$. Thus the covariance is unchanged, which is also what we would expect. This result easily generalizes (see Exercise 12) to any future increment of any size, which means that  the Brownian motion is independent of future increments. 
\close
\end{example}

\begin{example}
Let $x$ be a Brownian bridge. The covariance function is $\GG(s, t) = s\wedge t - st$. The formula for the conditional variance from Example~\ref{ex:cond_neumann} are now
\[
\begin{pmatrix}
s- s^2 & t - st \\
t - st & t - t^2
\end{pmatrix}
-\frac{\epsilon}{1-\epsilon}\begin{pmatrix}
s^2 &  st\\
st & t^2
\end{pmatrix}.
\]
As $\epsilon \rightarrow 0$, the covariance just remains that of the Brownian bridge.
\close
\end{example}

\section{Exercises}
\begin{itemize}
\item[1.] Consider the Poisson equation \eqref{eq:poisson} on the following three function spaces
\begin{align*}
\HH_1 &= \{f\in \mathcal{C}^2 \,:\, f(0)=f(1)=0\},\\[0.2em]
\HH_2 &= \{f\in \mathcal{C}^2 \,:\, f(0)=f'(1)=0\},\\[0.2em]
\HH_3 &= \{f\in \mathcal{C}^2 \,:\, f'(0)=f'(1)=0, {\textstyle \int_0^1f(t)\,\ud t = 0 }\}.
\end{align*}
Decide whether $-\partial^2_t f(t)$ can equal the constant function $1(t)=1$ for $f$ belonging to each of these sets. Explain the relation to examples~\ref{ex:indicator}~and~\ref{ex:bm_green}.
\item[2.] Show that 
\[
\HH = \{f\in \mathcal{C}^2([0,1])\,:\, {\textstyle \int_0^1f(t)\,\ud t = 0 }\}
\]
is a subspace of $\mathcal{C}^2([0,1])$.
\item[3.] Show that 
\[
\HH = \{f\in \mathcal{C}^2([0,1])\,:\, f(0)=f(1)\}
\]
is a subspace of $\mathcal{C}^2([0,1])$.
\item[4.] Show that 
\[
\HH_1 = \{f\in \mathcal{C}^2([0,1])\,:\, f(t) = f(1 - t) \text{ for all } t\in [0,1/2]\}
\]
and
\[
\HH_2 = \{f\in \mathcal{C}^2([0,1])\,:\, f(t) = -f(1 - t) \text{ for all } t\in [0,1/2]\}
\]
are subspaces of $\mathcal{C}^2([0,1])$.
\item[5.] Let $\HH$ be the linear subspace of functions in $\mathcal{C}^2([0, 1])$ for which $f(0) = f(1) = 0$ and $f^{(4)}(t) =0$. Identify the Green's function for the Laplace operator on $\HH$.
\item[6.] Let $\HH$ be the space of functions in $\mathcal{C}^2([0, 1])$ that obey homogeneous Dirichlet boundary conditions. Express this constraint in terms of Fourier coefficients. 
\item[7.] Let $x$ be a Brownian motion. Define a new process $\tilde{x}(t) = x(t) - t x(1) $. compute its mean and covariance function. What process is $\tilde{x}$?
\item[8.] Let $\tilde{x}$ be a backward Brownian motion, that is, $\tilde{x}(t) = x(1-t)$ for $t\in [0, 1]$ where $x$ is a standard Brownian motion. Determine the covariance function for $\tilde{x}$.
\item[9.] Let $x$ be the sum of independent forward and backward Brownian motions (see previous exercise). Determine the covariance function for $x$. Determine boundary conditions that this covariance function obeys.
\item[10.] Let $\tilde x$ be the sum of a forward and backward Brownian motion (see previous exercises), but this time assume that the backward Brownian motion arise from the forward Brownian motion, that is $\tilde x(t) = x(t) + x(1-t)$ where $x$ is the standard Brownian motion. Determine the covariance function for $x$. Determine boundary conditions that this covariance function obeys.

\item[11.] Argue why it must hold that the conditional expectation maximizes $\E[\theta\,|\, \by]$ the posterior $p(\theta\,|\, \by)$ when we assume that $\theta$ is a Gaussian process.
\item[12.] Show the general version of the result about conditioning on future increments for the Brownian motion mentioned in Example~\ref{ex:bm_cond_neumann}. You can either use the conditional distribution, or use the definiton of Brownian motion. 
\item[13.] Give an interpretation of the result in Example~\ref{ex:poly_motion} in terms of the Brownian motion.
\end{itemize}

\printbibliography

\typeout{get arXiv to do 4 passes: Label(s) may have changed. Rerun}
\end{document}